\documentclass[12pt,reqno]{amsart}

\usepackage[T1]{fontenc}
\usepackage{newpxtext}                     
\usepackage[varg,bigdelims]{newpxmath}    
\usepackage[mathscr]{eucal}                
\usepackage{bm}
\usepackage{microtype}

\usepackage[a4paper,textwidth=16cm,textheight=22cm]{geometry}
\usepackage{setspace}

\usepackage{amsmath,mathtools}    
\usepackage{graphicx}
\usepackage{xcolor}
\usepackage{ytableau}
\usepackage{stackengine}
\usepackage{enumitem}
\usepackage{tikz-cd}
\usetikzlibrary{arrows.meta}

\usepackage[backend=biber,style=alphabetic,maxbibnames=99,doi=false,url=false]{biblatex}

\usepackage[colorlinks,
            linkcolor={red!55!black},
            citecolor={green!40!black},
            urlcolor={blue!60!black}]{hyperref}
\usepackage[capitalise,noabbrev]{cleveref}

\tikzset{
  bulletepi/.style={
    {Circle[fill=black, scale=0.7]}-{Straight Barb[] Straight Barb[]}
  }
}

\theoremstyle{plain}
\newtheorem{theorem}{Theorem}[section]
\newtheorem{corollary}[theorem]{Corollary}

\newtheorem{lemma}[theorem]{Lemma}
\newtheorem{proposition}[theorem]{Proposition}

\newtheorem*{lemma*}{Technical Lemma}

\theoremstyle{definition}

\newtheorem{example}[theorem]{Example}
\newtheorem{rem}[theorem]{Remark}

\theoremstyle{plain}

\newcommand{\mcd}{\mathcal{D}_{\!\lambda\mu}^{\nu}}

\newcommand{\posreal}[1][+]{\Delta_{re}^{#1}}
\newcommand{\posimag}{\Delta_{im}^+}

\newcommand{\cato}{\mathcal{O}}
\newcommand{\catoint}{\mathcal{O}^{\mathrm{int}}}
\newcommand{\lie}{\mathfrak}

\newcommand{\lrw}{c_{\lambda\mu}^\nu(w)}
\newcommand{\lr}{c_{\lambda\mu}^\nu}

\DeclareMathOperator{\ann}{ann}

\DeclareMathOperator{\Hom}{Hom}

\begin{document}

\title{Splitting of multiplicity spaces into symmetric and alternating parts}
\author{Manika Gupta}
\address{The Institute of Mathematical Sciences, Chennai, India}
\address{Homi Bhabha National Institute, Training School Complex, Anushakti Nagar, Mumbai 400094, India}
\email{manikag@imsc.res.in}
\maketitle
\begin{abstract}
For a finite dimensional semisimple Lie algebra $\lie{g}$, the multiplicity $\lr$ of an irreducible $V(\nu),$ inside the tensor product $V(\lambda)\otimes V(\mu),$ of two highest weight irreducible representations can be realised as the dimension of the Kostant–Parthasarathy–Ranga Rao–Varadarajan multiplicity space inside $V(\lambda)$, $V(\mu)$ or $V(\nu)$. A result of Panyushev and Yakimova shows how to split the multiplicity space for $c_{\lambda\lambda}^{\nu}$ inside $V(\nu)$ into symmetric and alternating parts when $V(\lambda)$ is self-dual.
For finite-dimensional semisimple or symmetric Kac--Moody, we extend this to Kostant--Kumar modules $K(\lambda,w,\lambda)=\lie{Ug}(v_\lambda\otimes v_{w\lambda})$ when $w=w^{-1},$ and for all dominant weights $\lambda.$ We show that the multiplicities of an irreducible in the symmetric and alternating parts of $K(\lambda,w,\lambda)$ are given by the dimensions of the eigenspaces of the Tits lift $\widetilde w$ of $w$ on the multiplicity space corresponding to eigenvalues $\pm 1$. This recovers the result for the full tensor product by letting $w=w_0,$ the longest element of the Weyl group. In type A, we also discuss connections to a combinatorial interpretation due to Carr\'e and Leclerc in terms of domino tableaux. 
\end{abstract}

\section{Introduction}
Given a complex symmetrizable Kac--Moody Lie algebra $\lie g$, the problem of determining the multiplicity $\lr$ of an irreducible highest weight representation $V(\nu)$ inside the tensor product of two irreducible highest weight representations $V(\lambda)\otimes V(\mu)$ is well-studied. One of the ways to realize this multiplicity is through the dimension of the Kostant–Parthasarathy–Ranga Rao–Varadarajan (KPRV) multiplicity spaces:
$$
    \mathcal{D}_{\lambda\mu}^{\nu} \;=\; \bigl\{\,v\in V(\mu)_{\nu-\lambda}\;:\;
    e_i^{\,\langle\lambda,\alpha_i^\vee\rangle+1}\,v=0\ \text{ for all simple roots } \alpha_i\,\bigr\} \;\subseteq V(\mu)
$$

and when $\mathfrak g$ is finite dimensional, 

\[
    \mathcal{M}_{\lambda\mu}^{\nu} \;=\;\bigl\{\,v\in V(\nu)_{\lambda-\mu^*}\;:\;
    e_i^{\,\langle\mu^*,\alpha_i^\vee\rangle+1}\,v=0\ \text{ for all simple roots } \alpha_i\,\bigr\} \;\subseteq V(\nu).
\] 
That is, $\dim \mathcal{D}_{\lambda\mu}^{\nu}=\dim \mathcal{M}_{\lambda\mu}^{\nu}=\lr.$

For $\lambda=\mu$, we have $V(\lambda)\otimes V(\lambda)= \text{Sym}^2 V(\lambda)\oplus \text{Alt}^2V(\lambda)$ and so, $c_{\lambda\lambda}^\nu=c_{\lambda\lambda}^{\nu}(+) + c_{\lambda\lambda}^{\nu}(-),$ where $c_{\lambda\lambda}^{\nu}(+)$ and $c_{\lambda\lambda}^{\nu}(-)$ are the multiplicities of $V(\nu)$ in $\text{Sym}^2 V(\lambda)$ and $\text{Alt}^2V(\lambda)$ respectively. A result of Panyushev and Yakimova shows how to split the multiplicity space for $c_{\lambda\lambda}^{\nu}$ into symmetric and alternating parts in the finite dimensional case where $\lambda=\lambda^*,$ in terms of the Tits lift of $w_0$, the longest element of the Weyl group and the Frobenius-Schur indicator of $V(\lambda).$

Kostant--Kumar modules are submodules of the tensor product, indexed by elements $w$ of the Weyl group of $\lie g$:
\[
    K(\lambda,w,\mu)\;:=\;\lie{Ug}\,(v_\lambda\otimes v_{w\mu})
    \;\subseteq\; V(\lambda)\otimes V(\mu),
\]
where $v_\lambda$ and $v_{w\mu}$ are nonzero vectors of weights $\lambda$
and $w\mu$ in $V(\lambda)$ and $V(\mu)$ respectively. 
If we restrict to the case when $\lambda=\mu$ and $w=w^{-1},$ $$K(\lambda, w, \lambda)= (K(\lambda, w, \lambda)\cap \text{Sym}^2V(\lambda)) \oplus (K(\lambda, w, \lambda)\cap \text{Alt}^2V(\lambda)).$$ 

For a symmetrizable Kac--Moody Lie algebra, we define an involution on $\mathcal{D}_{\lambda\lambda}^{\nu}$ such that the dimension of its $\pm 1$ eigenspaces count the multiplicities $c_{\lambda\lambda}^{\nu}(+)$ and $ c_{\lambda\lambda}^{\nu}(-).$

For finite-dimensional semisimple or symmetric Kac--Moody, using a generators and relations description of Kostant--Kumar modules, we can define a KPRV-like multiplicity space inside $V
(\nu)$, $$\mathcal{M}_{\lambda\lambda}^{\nu}(w)=\{v \in V(\nu)_{\lambda+w\lambda}: e_\alpha^{p_\alpha(w\lambda)+1} v =0 \text{ for } \alpha \in \posreal \text{ and } {\lie g}_{\beta} v =0 \text{ for } \beta \in \posimag\}.$$ This multiplicity space is left invariant by the action of the Tits lift $\widetilde{w}$ of $w.$ 

The main result of this paper shows how this multiplicity space splits into symmetric and alternating parts. The  multiplicities of an irreducible $V(\nu)$ in $(K(\lambda, w, \lambda)\cap \text{Sym}^2V(\lambda))$ and $(K(\lambda, w, \lambda)\cap \text{Alt}^2V(\lambda))\;$ are given by the dimensions of the eigenspaces corresponding to eigenvalues $\varepsilon_{w\lambda}$ and $-\varepsilon_{w\lambda}$ of the operator $\widetilde {w}$ on the space $\mathcal{M}_{\lambda\lambda}^{\nu}(w),$ where $\widetilde{w}^2(v_{\lambda})=\varepsilon_{w\lambda}v_\lambda,$ and $\varepsilon_{w\lambda}$ is a scalar with value either $1$ or $-1.$ 

Finally, in type $A$, we discuss combinatorial interpretations of the dimensions of the multiplicity spaces using a result due to Carré and Leclerc.\\

\noindent\textbf{AI tools use disclosure:} LLMs were used to find typos, errors, and check references. The mathematical arguments and proofs in this paper were entirely human-produced, with no involvement of AI tools at any stage.

\section{Notation and preliminaries}
Let $\lie g$ be a symmetrizable Kac--Moody Lie algebra over $\mathbb{C}$ with triangular decomposition $\mathfrak{g}=\mathfrak{n}^{-}\oplus \mathfrak{h} \oplus\mathfrak{n}^{+}$. Let $\mathfrak{b}=\mathfrak h\oplus \mathfrak{n}^+$ be its Borel subalgebra. 
If the generalised Cartan matrix of $\lie g$ is symmetric, we call $\lie g$ a symmetric Kac--Moody Lie algebra. Let $\alpha_i, \,i=1,2,\ldots,n$ be the simple roots of $\mathfrak{g}$. Let $\Delta^+$ and $\Delta^-$ be the sets of positive and negative roots respectively and $\posreal, \,\posimag$ denote the sets of positive real and imaginary roots. Let $P$ and $P^+$ be the sets of weights and dominant weights respectively. Let $\rho^\vee\in \lie h$ denote an element that satisfies $\langle\alpha,\rho^\vee\rangle=1$ for all simple roots $\alpha$ of $\lie g$.

For each root $\beta$, let ${\lie g}_\beta$ denote the corresponding root space of ${\lie g}$. Given $\alpha \in \posreal$, we let $e_\alpha, f_\alpha$ denote non-zero elements of ${\lie g}_\alpha$ and ${\lie g}_{-\alpha}$ respectively. Let $e_i=e_{\alpha_i}$ and $f_i=f_{\alpha_i}$ for simple roots $\alpha_i.$ Let $W$ be the Weyl group of $\lie g,$ generated by simple reflections $s_1,s_2,\ldots,s_n.$

Let $\catoint$ denote the category of integrable $\lie g$-modules in category $\cato$. For a dominant integral weight $\lambda$, let $V(\lambda)$ denote the
irreducible $\lie g$-module of highest weight $\lambda$.

We denote the multiplicity of an irreducible $V(\nu)$ in  $V(\lambda)\otimes V(\mu)$ and $K(\lambda,w,\mu)$ as $\lr$ and $\lr(w)$ respectively. 

For $\gamma \in P,$ the height, denoted by $ht(\gamma)$ is the sum of the coefficients when $\gamma$ is expressed in terms of the simple roots. Hence, $ht(\gamma)= \langle \gamma,\rho^\vee\rangle.$

 An integrable highest weight module $V(\lambda)$ of $\lie g$ admits a positive definite Hermitian form $\langle \cdot, \cdot \rangle$ (the {\em Shapovalov form}) satisfying the $\omega_0$-contravariance property:
 \[ \langle Xv, v' \rangle = \langle v, \omega_0(X) v' \rangle \]
    for all $v, v' \in V$ and $X \in \lie g$.
This form on $V(\lambda)$ is unique up to scaling. If $V = V(\lambda_1) \otimes V(\lambda_2)$ for some $\lambda_1, \lambda_2 \in P^+$, then we consider the following Shapovalov form: $V(\lambda_1)$ and $V(\lambda_2)$:
\[\langle v_1\otimes v_2, v_1'\otimes v_2'\rangle = \langle v_1,v_1' \rangle \langle v_2,v_2' \rangle\]
for $v_i, v_i' \in V(\lambda_i)$ for $i=1, 2$.

\section{KPRV multiplicity space and Tits lifts}

In this section, we recall Kostant-Parthasarathy-Ranga Rao-Varadarajan multiplicity spaces and some properties of Tits lifts of elements of the Weyl group. 
\subsection{} When $\lie g$ is a symmetrizable Kac--Moody Lie algebra, the multiplicity of an irreducible $V(\nu)$ in a tensor product can be realised as the dimension of the space:
\begin{equation}\label{eq:prv1mult}
    \mathcal{D}_{\lambda\mu}^{\nu} \;=\; \bigl\{\,v\in V(\mu)_{\nu-\lambda}\;:\;
    e_i^{\,\langle\lambda,\alpha_i^\vee\rangle+1}\,v=0\ \text{ for all simple roots } \alpha_i\,\bigr\} \;\subseteq V(\mu).
\end{equation}
This is due to  Parthasarathy-Ranga Rao-Varadarajan \cite{PRV} and Kostant \cite{kostant-tensmult} in the finite-dimensional case and was extended to all symmetrizable Kac--Moody algebras by Jeralds and Kumar \cite{jeralds-kumar}.

In the finite dimensional case, using the hom-tensor duality, we get $c_{\lambda\mu}^{\nu}=c_{\mu^*\nu}^{\lambda},$ where $\lambda^*=-w_0\lambda$ is the dual of $\lambda.$
Defining \begin{equation}\label{eq:mult}
    \mathcal{M}_{\lambda\mu}^{\nu}= \bigl\{\,v\in V(\nu)_{\lambda-\mu^*}\;:\;
    e_i^{\,\langle\mu^*,\alpha_i^\vee\rangle+1}\,v=0\ \text{ for all simple roots } \alpha_i\,\bigr\} \subseteq V(\nu),
\end{equation}
we obtain $\lr=\dim\mathcal{M}_{\lambda\mu}^{\nu}.$
For finite dimensional semisimple or symmetric Kac--Moody $\lie g$ and $w\in W$, we will define an analogous multiplicity space $\mathcal{M}_{\lambda\mu}^{\nu}(w)$ later.

\subsection{}
Let $\lie g$ be a symmetrizable Kac--Moody Lie algebra and $\phi: \lie g\to \text{End}\;V$ be a representation of $\lie g$ in $\catoint$. For a simple reflection $s_i$ in $W,$ its Tits lift on $V$ is defined as $$\widetilde{s_i}= \text{exp }\phi(e_i)\; \text{exp }(-\phi(f_i))\; \text{exp }\phi(e_i)\in GL(V).$$ Given an element $w\in W,$ with reduced expression $s_{i_1}\ldots s_{i_k},$ its Tits lift $\widetilde{w}$ is defined as $\widetilde{s}_{i_1}\ldots\widetilde{s}_{i_k}.$ This is well-defined because the lifts for simple reflections satisfy braid relations, see \cite{T1}, \cite{S1} and Chapter 6 of \cite{Kum3}.

For a simple reflection $s_i,$ the square of its lift acts as \begin{equation}\label{eq:simplelift}
    (\widetilde{s_i})^2v=(-1)^{\langle \mu,\alpha_i^\vee\rangle}v 
\end{equation} where the weight of $v$ is $\mu.$

For $\gamma\in P,$ define $\varepsilon_{w\gamma}=(-1)^{ht(\gamma-w\gamma)}.$

Then, more generally, we have the following well-known result. We sketch the proof for completeness. 
\begin{proposition}\label{proptitslift}
    Let $\lie g$ be a symmetrizable Kac--Moody Lie algebra and $V$ be a representation of $\lie g$ in $\catoint$. Let $w\in W$ be such that $w^2=\text{id}.$ Then, $
\widetilde w^2\,v=\varepsilon_{w\mu}v,$ where the weight of $v$ is $\mu$.
\end{proposition}
\begin{proof}
    Let $w\in W$ have reduced expression $s_{i_1}\ldots s_{i_k}$. Define $t_{\alpha_i^\vee}=\text{exp}(\pi i \alpha_i^\vee) \in GL(V).$ From Equation ~\eqref{eq:simplelift}, we have that  $(\widetilde{s_i})^2=t_{\alpha_i^\vee}.$ \\
    $\widetilde{w}^2=\widetilde{s_{i_1}}\ldots\widetilde{s_{i_k}}\widetilde{s_{i_k}}\ldots\widetilde{s_{i_1}}= \widetilde{s_{i_1}}\ldots\widetilde{s_{i_{k-1}}}t_{\alpha_{i_k}^\vee}\widetilde{s_{i_{k-1}}}\ldots\widetilde{s_{i_1}}.$ Also, $\widetilde{s_i}\,t_\lambda\,\widetilde{s_i}^{-1} =t_{s_i\lambda}, $ see \cite{S1}. So, we have
$\widetilde{s_i}\,t_{\lambda}\,\widetilde{s_i}
=\bigl(\widetilde{s_i}\,t_\lambda\,\widetilde{s_i}^{-1}\bigr)\,\widetilde{s_i}^{\,2}
=t_{s_i\lambda}\,t_{\alpha_i^\vee}
=t_{s_i\lambda+\alpha_i^\vee}.$ 

\[
\widetilde w^2=\widetilde{s_{i_1}}\cdots\widetilde{s_{i_{k-1}}}\;t_{\alpha_{i_k}^\vee}\;\widetilde{s_{i_{k-1}}}\cdots\widetilde{s_{i_1}}.
\]
Apply the identity to the innermost $\widetilde{s_{i_{k-1}}}\,t_{\alpha_{i_k}^\vee}\,\widetilde{s_{i_{k-1}}}$ and continue outwards. Note that each step applies $s_{i_m}$ to everything already in the subscript and adds a new term $\alpha_{i_m}^\vee$. After $k-1$ steps, we have
\[
\widetilde w^2=t_{\lambda},\qquad
\lambda=\alpha_{i_1}^\vee+s_{i_1}\alpha_{i_2}^\vee+s_{i_1}s_{i_2}\alpha_{i_3}^\vee+\cdots+s_{i_1}\cdots s_{i_{k-1}}\alpha_{i_k}^\vee
=\sum_{j=1}^k\beta_j^\vee,
\]
where $\beta_j=s_{i_1}\cdots s_{i_{j-1}}\alpha_{i_j}.$

Since the word is reduced,
$\{\beta_1,\dots,\beta_k\}=\{\beta>0:\ w^{-1}\beta<0\}$, each root occurring once. As $w=w^{-1}$, this set equals $\{\beta>0:\ w\beta<0\}$.

Since $t_\lambda=\exp(\pi i\lambda)$ acts on $V_\mu$ by $(-1)^{\langle\mu,\lambda\rangle}$, for $v$ of weight $\mu$ we get
$$
\widetilde w^2\,v=(-1)^{\sum_{\beta>0,\ w\beta<0}\langle\mu,\beta^\vee\rangle}\,v=(-1)^{\langle \mu, \rho^\vee-w\rho^\vee\rangle}\,v=(-1)^{\langle \mu-w\mu, \rho^\vee\rangle}\,v.
$$
\end{proof}
\begin{corollary}\label{cor:eigen}
    Let $\lie g$ be a symmetrizable Kac--Moody Lie algebra and $V$ be a representation of $\lie g$ in $\catoint$. Let $w\in W$ be such that $w^2=\text{id}.$ If $\mu\in P$ satisfies $w\mu=\mu,$ then $\widetilde{w}^2$ acts as $\text{id}$ on $V_\mu,$ the $\mu$ weight space of $V.$ In particular, the eigenvalues of $\widetilde{w}$ on $V_\mu$ belong to $\{\pm1\}.$
\end{corollary}

\subsection{}
Let $\lie g$ be  finite dimensional and $V(\lambda)$ be self-dual. We have a non-degenerate $\lie g$-invariant bilinear form $\Phi_{\lambda}$ on $V({\lambda})$. This form is either symmetric or skew-symmetric, that is,   $\Phi_{\lambda}(x,y)=\epsilon_\lambda\Phi_{\lambda}(y,x)$ for all
$x,y\in V({\lambda})$, where $\epsilon_\lambda\in\{\pm1\}.$ The sign $\epsilon_\lambda$ is called the Frobenius-Schur indicator of $V(\lambda).$
We have the following well-known relationship between the Tits lift of the longest element $w_0\in W$ and the Frobenius-Schur indicator of $V(\lambda)$.

\begin{proposition}
    Let $\lie g$ be a finite dimensional semisimple Lie algebra, $\lambda\in P^+$ and let $V(\lambda)$ be self-dual, that is,  $V(\lambda)\cong V(\lambda)^*.$ Then the Frobenius-Schur indicator of $V(\lambda)$ is equal to the eigenvalue of $\widetilde{w_0}^2$ on the weight space $\mathbb{C}v_\lambda,$ which is equal to $\varepsilon_{w_0\lambda}.$
\end{proposition}
\begin{proof}
    Define $\widetilde{w_0}v_\lambda=:v_{w_0\lambda}$. By Proposition ~\ref{proptitslift}, $\widetilde{w_0}^2v_\lambda=\varepsilon_{w_0\lambda}v_\lambda.$ Since $\Phi_\lambda$ is invariant under $\text {exp} \; \lie g,$ we have, $$\Phi_\lambda(v_\lambda,v_{w_0\lambda})=\Phi_\lambda(\widetilde{w}_0v_\lambda,\widetilde{w}_0v_{w_0\lambda})= \varepsilon_{w_0\lambda} \Phi_\lambda(v_{w_0\lambda}, v_\lambda).$$ The non-degeneracy of $\Phi_\lambda$ ensures that $\Phi_\lambda(v_\lambda,v_{w_0\lambda})$ is non-zero.
    Therefore, the Frobenius-Schur indicator, $\epsilon_\lambda$ is equal to $\varepsilon_{w_0\lambda}.$
    
\end{proof}

\section{Kostant--Kumar modules}
\subsection{}
Let $\lie g$ be a symmetrizable Kac--Moody algebra and let $\lambda, \mu \in P^+$. Kostant--Kumar modules are defined as $\mathfrak{Ug}$-submodules of the tensor product $V=V(\lambda)\otimes V(\mu)$, generated by a single vector of the form $v_\lambda\otimes v_{w\mu},$ where $w\in W,$ the Weyl group of the Lie algebra, and $v_\lambda$ is a non-zero vector of weight $\lambda$ in $V(\lambda)$ and $v_{w\mu}$ is a non-zero vector of weight $w\mu$ in $V(\mu).$ They are denoted by \[K(\lambda,w,\mu)=\mathfrak{Ug}(v_{\lambda}\otimes v_{w\mu}).\]
When $\lie g$ is a finite dimensional semisimple Lie algebra, setting $w=w_0,$ the longest element of the Weyl group, the module $K(\lambda, w_0, \mu)$ becomes the entire tensor product. 

We use the following generators and relations description of these modules given in Theorem 4.2 of \cite{MRV} to define a multiplicity space inside $V(\nu)$.
\begin{theorem}\label{thm:gen-rel}
    Let $\mathfrak{g}$ be a finite dimensional semisimple Lie algebra or a symmetric Kac--Moody Lie algebra. Consider the Kostant--Kumar module $K(\lambda,w,\mu)=\mathfrak{Ug} (v_\lambda\otimes v_{w\mu})$. Given $\gamma \in P$ and $\alpha \in \posreal$, define 
\begin{align*}
    p_\alpha(\gamma) &:= \max \{0, -\langle \gamma, \alpha^\vee\rangle\}\\
    q_\alpha(\gamma) &:= \max \{0, \langle \gamma, \alpha^\vee\rangle\} 
\end{align*} 
Let $J$ be the left ideal of $\mathfrak{Ug}$ generated by the elements:  
    \begin{enumerate}
        \item $e_\alpha^{\,p_\alpha (w\mu)+1}, \;\; \alpha \in \posreal$.
        \smallskip\item $f_\alpha^{\,\langle \lambda, \alpha^\vee \rangle + q_\alpha (w\mu)+1}, \;\; \alpha \in \posreal$.
        \smallskip\item ${\lie g}_\beta, \;\; \beta \in \posimag$.
        \smallskip\item $ h - \langle \lambda + w\mu, h \rangle 1, \;\; h \in \mathfrak{h}$.
    \end{enumerate}
    Then $J=\ann_{\,\mathfrak{Ug}} (v_\lambda \otimes v_{w\mu})$, and $K(\lambda, w, \mu) \cong \mathfrak{Ug}/J$. \qed
\end{theorem}
Define $$\mathcal{M}_{\lambda\mu}^\nu(w)= \{v \in V(\nu)_{\lambda+w\mu}: e_\alpha^{p_\alpha(w\mu)+1} v =0 \text{ for } \alpha \in \posreal \text{ and } {\lie g}_{\beta} v =0 \text{ for } \beta \in \posimag\}.$$

\begin{corollary}
    Under the hypotheses of the previous theorem, the map \begin{equation} \Hom_{\lie g} (K(\lambda,w,\mu), V(\nu))\to \mathcal{M}_{\lambda\mu}^{\nu}(w) \end{equation}\label{map}
    where $f\mapsto f(v_\lambda\otimes v_{w\mu})$ is an isomorphism. Therefore, we have 
    $\lrw = \dim \mathcal{M}_{\lambda\mu}^\nu(w).$
\end{corollary}
\begin{proof}
    A map $f\in \Hom_{\lie g} (K(\lambda,w,\mu), V(\nu))$ is determined by its value on the generator $v_\lambda\otimes v_{w\mu}.$ Since $f: \lie{Ug}/J\to V(\nu),$ $J\subseteq \ann_{\lie{Ug}}f(v_\lambda\otimes v_{w\mu}).$ Note that by Lemma ~\ref{lem:sl2string}, for any vector $v\in \mathcal{M}_{\lambda\mu}^{\nu}(w),$ $f_\alpha^{\,\langle \lambda, \alpha^\vee \rangle + q_\alpha (w\mu)+1}v=0, \;\forall\; \alpha \in \posreal$. Conversely, given any element $v'$ in $\mathcal{M}_{\lambda\mu}^{\nu},$ we can define $f\in \Hom_{\lie g} (K(\lambda,w,\mu), V(\nu))$ by setting $f(v_\lambda\otimes v_{w\mu})=v'.$
\end{proof}

\begin{rem}\label{rem:mult}
    When $\lie g$ is finite dimensional, $\mathcal{M}_{\lambda\mu}^\nu(w_0)=\mathcal{M}_{\lambda\mu}^\nu$ as defined in Equation ~\eqref{eq:mult} (see \cite{MRV}).
\end{rem}

\subsection{}

Let $\tau: V(\lambda)\otimes V(\mu)\to V(\mu)\otimes V(\lambda)$ be  the swap operator, that is $\tau(v\otimes v')=v'\otimes v,$ extended linearly. Recall (see \cite{KRV} and \cite{KRV2}) that $\tau(K(\lambda,w,\mu))=K(\mu,w^{-1},\lambda).$ Hence, for $\lambda=\mu,$ and $w$ an involution, $K(\lambda,w,\lambda)\subseteq V(\lambda)\otimes V(\lambda)$ is $\tau-$invariant.

Now, since, $$V(\lambda)\otimes V(\lambda)=\text{Sym}^2V(\lambda)\oplus \text{Alt}^2V(\lambda),$$ we obtain $$K(\lambda, w, \lambda)= (K(\lambda, w, \lambda)\cap \text{Sym}^2V(\lambda)) \oplus (K(\lambda, w, \lambda)\cap \text{Alt}^2V(\lambda)).$$ 
Let $K^{+}(\lambda,w, \lambda)=K(\lambda, w, \lambda)\cap \text{Sym}^2V(\lambda)$ and $K^{-}(\lambda,w, \lambda)=K(\lambda, w, \lambda)\cap \text{Alt}^2V(\lambda).$

\subsection{}
We now state two well-known lemmas (see \cite{Hum1} and Section 3.8 of \cite{Kac}) which we will use to prove the invariance of the multiplicity space $\mathcal{M}_{\lambda\lambda}^\nu(w)$ with respect to the action of $\widetilde{w}.$

\begin{lemma}\label{lem:sl2string}
Let $\mathfrak{g}$ be a symmetrizable Kac--Moody algebra, $V \in \catoint$ be a $\mathfrak{g}$-module and $\alpha \in \posreal$. Let $v \in V$ be a weight vector of weight $\gamma$ and let $k \geq 0$. Then $e_\alpha^k(v) =0$ if and only if $f_\alpha^{\langle \gamma, \alpha^\vee \rangle + k}\, v =0$.
\end{lemma}

\begin{lemma}
Let $\mathfrak{g}$ be a symmetrizable Kac--Moody algebra, $\phi: \lie g \to \text{End}\;V$ be a $\mathfrak{g}$-module in $\catoint$. Consider the operator $\widetilde{s_i}$ on $V,$ and let $x \in \lie g.$ Then, $$\widetilde{s_i}\phi(x)\widetilde{s_i}^{-1}=\phi(s_i^{ad}x),$$ where $s_i^{ad}= exp(ad(e_i))exp(ad(-f_i))exp(ad(e_i))\in \text{End}\; \lie g.$
\end{lemma}

\begin{corollary}\label{cor:twistkill}
    Let $\mathfrak{g}$ be a symmetrizable Kac--Moody algebra, $w\in W$ and $V$ be a $\mathfrak{g}$-module in $\catoint$. Then $e_\alpha^k(\widetilde{w}v)=0\iff e_{w^{-1}\alpha}^k(v)=0.$ (Here $e_{-\alpha}=f_{\alpha}.$)
\end{corollary}

\begin{proposition}\label{prop:invar}
    Let $\lie g$ be a symmetrizable Kac--Moody Lie algebra. Let $\lambda$ and $\nu$ be dominant weights and $w$ be an element of the Weyl group such that $w=w^{-1}$. The space $$ \mathcal{M}_{\lambda\lambda}^\nu(w)=\bigl\{v \in V(\nu)_{\lambda+w\lambda}: e_\alpha^{p_\alpha(w\lambda)+1} v =0 \text{ for } \alpha \in \posreal \text{ and } {\lie g}_{\beta} v =0 \text{ for } \beta \in \posimag\bigr\}$$ is invariant under the action of $\widetilde{w}.$ 
\end{proposition}
\begin{proof}
    Let $v\neq 0$ be an element of $\mathcal{M}_{\lambda\lambda}^\nu(w).$
    
    $\widetilde{w}$ sends the weight space corresponding to $\lambda+w\lambda$ to $w(\lambda+w\lambda)=w\lambda+\lambda.$ 
    
    If $w\alpha>0,$ by Corollary ~\ref{cor:twistkill}, $e_{w^{-1}\alpha}^mv=e_{w\alpha}^mv=0 \implies  e_{\alpha}^m(\widetilde{w}v)=0.$ Now, $p_{w\alpha}(w\lambda)+1= \text{max}\{0, -\langle w\lambda, w\alpha^\vee\rangle\}+1=1,$ since $\lambda$ is dominant. So $e_\alpha(\widetilde{w}v)=0\implies e_\alpha^{p_\alpha(w\lambda)+1}(\widetilde{w}v)=0.$
    
    If $w\alpha< 0,$ consider $p_\alpha(w\lambda)+1= \max\{0,-\langle w\lambda,\alpha^\vee\rangle\}+1=-\langle\lambda, w\alpha^\vee\rangle+1,$ since $w\alpha$ is a negative root and $\lambda$ is dominant. By Corollary ~\ref{cor:twistkill}, to show $e_{\alpha}^{-\langle\lambda,w\alpha^\vee\rangle+1}\widetilde{w}v=0,$ we need to show that $f_{-w\alpha}^{-\langle\lambda,w\alpha^\vee\rangle+1}v=0.$ Consider the $-w\alpha$ string through $v.$ We know that $e_{-w\alpha}^kv=0$ where $k=p_{-w\alpha}(w\lambda)+1=\langle\lambda,\alpha^\vee\rangle+1.$ Note that $v$ has weight $\lambda+w\lambda,$ and let $n=\langle\lambda+w\lambda,-w\alpha^\vee\rangle.$ Let $k'$ be such that $f^{k'+1}_{-w\alpha}v=0,$ and $f^{k'}_{-w\alpha}v\neq 0.$ Then, $n\ge k'-\langle\lambda,\alpha^\vee\rangle,$ by Lemma ~\ref{lem:sl2string}, which gives us $-\langle\lambda, w\alpha^\vee\rangle\ge k',$ so we are done.
    
    Clearly, the action of $\widetilde{w}$ preserves the condition ${\lie g}_{\beta} v =0 \text{ for } \beta \in \posimag$ since the Weyl group preserves the set of positive imaginary roots.
\end{proof}

\section{Splitting of multiplicity spaces}

Considering $V(\lambda)\otimes V(\lambda)$ over a symmetrizable Kac--Moody Lie algebra $\lie g$ prompts the natural question of determining the multiplicity of an irreducible in $\text{Sym}^2 V(\lambda)$ and $\text{Alt}^2V(\lambda).$ If $c_{\lambda\lambda}^{\nu}(+)$ and $c_{\lambda\lambda}^{\nu}(-)$ are the multiplicities of $V(\nu)$ in $\text{Sym}^2 V(\lambda)$ and $\text{Alt}^2V(\lambda)$ respectively, we have $c_{\lambda\lambda}^\nu=c_{\lambda\lambda}^{\nu}(+) + c_{\lambda\lambda}^{\nu}(-).$ This can be realised at the level of multiplicity spaces as well. 

\subsection{}
We have the following result about the multiplicity space $\mcd$ defined in Equation ~\eqref{eq:prv1mult}.
(See \cite{GE} and \cite{MRV}.)
\begin{proposition}\label{prop:pi0-image}
    Let $\mathfrak{g}$ be a symmetrizable Kac--Moody algebra and $\lambda, \mu, \nu \in P^+$. 
    Let $v_\lambda$ denote a nonzero highest weight vector of $V(\lambda)$. Let $\pi_0: V(\lambda)\otimes V(\mu)\twoheadrightarrow v_\lambda\otimes V(\mu)$ be the orthogonal projection with respect to the product Shapovalov form. Let $\left(V(\lambda)\otimes V(\mu)\right)_{\nu}^+$  be the set of highest weight vectors of weight $\nu$ in the tensor product. Then, $\pi_0(\left(V(\lambda)\otimes V(\mu)\right)_\nu^+)= v_{\lambda}\otimes \mcd.$
\end{proposition}
 
Define $\sigma: v_\lambda\otimes \mathcal{D}_{\lambda\lambda}^\nu\to v_\lambda\otimes \mathcal{D}_{\lambda\lambda}^\nu$ as $\sigma:=\pi_0\tau\pi_0^{-1}$ where $\tau$ is the swap operator. Let $\mathcal{D}_{\lambda\lambda}^{\nu}(\sigma,\pm1)$ denote the eigenspaces of the operator $\sigma$ on $v_\lambda\otimes \mathcal{D}_{\lambda\lambda}^{\nu}$ corresponding to the eigenvalues $\pm 1$ respectively. These eigenspaces count the multiplicities of symmetric and alternating parts as follows.

\begin{theorem}\label{main}
    Let $\lie g$ be a symmetrizable Kac--Moody Lie algebra. Let $\lambda$ and $\nu$ be dominant weights. Then,
    \[
c^{\nu}_{\lambda\lambda}(+)=\dim \mathcal{D}_{\lambda\lambda}^{\nu}(\sigma,+1)
\qquad
c^{\nu}_{\lambda\lambda}(-)=\dim\mathcal{D}_{\lambda\lambda}^{\nu}(\sigma,-1).
\] 
\end{theorem}
\begin{proof}
    If $v\in (V(\lambda)\otimes V(\lambda))_\nu^+ \;\bigcap \;\text{Sym}^2V(\lambda),$ then, $\tau(v)=v.$ Hence $\pi_0(v)$ is an eigenvector of $\sigma$ corresponding to eigenvalue $+1.$ Conversely, if $v'\in v_\lambda\otimes \mathcal{D}_{\lambda\lambda}^{\nu}$ is an eigenvector corresponding to eigenvalue $+1,$ then $\pi_0^{-1}(v')\in \text{Sym}^2V(\lambda).$ Hence, $c^{\nu}_{\lambda\lambda}(+)=\dim \mathcal{D}_{\lambda\lambda}^{\nu}(\sigma,+1).$ A similar argument follows for the alternating case.
\end{proof}
\begin{rem}
    We can define an involution on the space $\mathcal{D}_{\lambda\lambda}^{\nu}$ instead of $v_\lambda\otimes \mathcal{D}_{\lambda\lambda}^{\nu}$ by `ignoring' the first component.
\end{rem}

\subsection{} Let $w\in W$ be an involution. Consider the Kostant--Kumar module $K(\lambda,w,\lambda)=\mathfrak{Ug} (v_\lambda\otimes v_{w\lambda}).$ 
Let the multiplicity of an irreducible $V(\nu)$ in $K^+(\lambda,w,\lambda)$ and $K^-(\lambda,w,\lambda)$ be denoted by $c_{\lambda\lambda}^{\nu}(w,+)$ and $c_{\lambda\lambda}^{\nu}(w,-)$ respectively.
Then, we also have the following result for Kostant--Kumar modules.
\begin{theorem}\label{mainw}
    Let $\lie g$ be a symmetrizable Kac--Moody Lie algebra. Let $\lambda$ and $\nu$ be dominant weights and $w\in W$ be such that $w^2=\text{id}. $ Let $\pi_w:v_\lambda\otimes V(\lambda)\to v_\lambda\otimes V_w(\lambda)$ be the orthogonal projection onto the Demazure module with respect to the Shapovalov form. Then,
    \[
c^{\nu}_{\lambda\lambda}(w,+)=\dim \pi_w(\mathcal{D}_{\lambda\lambda}^{\nu}(\sigma,+1))
\qquad
c^{\nu}_{\lambda\lambda}(w,-)=\dim\pi_w(\mathcal{D}_{\lambda\lambda}^{\nu}(\sigma,-1)).
\] 
\end{theorem}
\begin{proof}
    Let $u\in (V(\lambda)\otimes V(\lambda))_{\nu}^{+}$ be such that $\pi_w\pi_0(u)=0.$ Let $z\in V_w(\lambda)$. Since $\pi_0$ and $\pi_w$ are orthogonal projections,
    $\langle u, v_\lambda\otimes z\rangle=\langle \pi_0 u, v_\lambda\otimes z\rangle=\langle \pi_w\pi_0 u, v_\lambda\otimes z\rangle=0.$ Further, $K(\lambda,w,\lambda)=\lie{Un}^-(v_\lambda\otimes V_w(\lambda)).$ Let $y \in \lie n^-,$ then $\langle  u, y(v_\lambda\otimes z)\rangle=\langle xu,v_\lambda\otimes z\rangle$ where $x\in \lie n^+.$ Since $u$ is a highest weight vector, $x.u=0,$ and $\langle xu,v_\lambda\otimes z\rangle=0.$ Therefore, $\langle u, K(\lambda,w,\lambda)\rangle=0.$
    
    Let $K(\lambda,w,\lambda)_{\nu}^{+}$ be the space of highest weight vectors of weight $\nu$ in $K(\lambda,w,\lambda).$
    Let $v\in \pi_w( \mathcal{D}_{\lambda\lambda}^{\nu}(\sigma,+1))$ Theorem 3.8 of \cite{MRV} shows that $K(\lambda,w,\lambda)_{\nu}^{+}$ is mapped isomorphically via $\pi_w\pi_0$ to $ \pi_w(v_\lambda\otimes \mathcal{D}_{\lambda\lambda}^{\nu}).$ 
    Also, the previous theorem shows $\pi_0((V(\lambda)\otimes V(\lambda))_{\nu}^{+}\cap \text{Sym}^2V(\lambda)) = \mathcal{D}_{\lambda\lambda}^{\nu}(\sigma,+1).$ 
    Hence, $v=\pi_w\pi_0(v_1)=\pi_w\pi_0(v_2),$ where $v_1$ and $v_2$ are elements of $K(\lambda,w,\lambda)_{\nu}^{+}$ and  $\text{Sym}^2V(\lambda)$ respectively. Now, $\pi_w\pi_0(v_1-v_2)=0,$ and $v_1-\tau(v_1)=v_1+v_2- \tau(v_1)-v_2=(v_1-v_2)-\tau(v_1-v_2).$ Therefore, $v_1-\tau(v_1)\in K(\lambda,w,\lambda)$ and $\langle\tau(v_1-v_2), K(\lambda,w,\lambda)\rangle=\langle v_1-v_2, K(\lambda,w,\lambda)\rangle=0$ because the product Shapovalov form on a tensor product is preserved by $\tau.$ Since the form is positive definite, $v_1-\tau(v_1)=0.$ Hence, $v_1\in \text{Sym}^2V(\lambda).$ Hence $v\in \pi_w\pi_0(K(\lambda,w,\lambda)^+_{\nu}\cap\text{Sym}^2V(\lambda)).$
    
    Doing a similar computation with $v'\in \pi_w( \mathcal{D}_{\lambda\lambda}^{\nu}(\sigma,-1))$ gives $v'\in \pi_w\pi_0(K(\lambda,w,\lambda)^+_{\nu}\cap\text{Alt}^2V(\lambda)). $

\end{proof}

\subsection{}
Consider the case where $\lie g$ is a finite-dimensional semisimple Lie algebra and $\lambda=-w_0\lambda,$ that is, \textit{$\lambda$ is self-dual.}
In this case, it follows that (see Remark ~\ref{rem:mult})  $$\mathcal{M}_{\lambda\lambda}^{\nu}=\mathcal{M}_{\lambda\lambda}^{\nu}(w_0)=\{\,v\in V(\nu)_{0}\;:\;
    e_i^{\,\langle\lambda,\alpha_i^\vee\rangle+1}\,v=0\ \text{ for all simple roots } \alpha_i\,\}$$ and 
$
c^{\nu}_{\lambda\lambda}=\dim\mathcal{M}_{\lambda\lambda}^{\nu}(w_0).
$

It follows from Proposition ~\ref{prop:invar} that  $\mathcal{M}_{\lambda\lambda}^{\nu}(w_0)$ is
$\widetilde{w_0}$-stable. By Corollary ~\ref{cor:eigen}, the eigenvalues of $\widetilde{w_0}$ on $V(\nu)_0$ belong to $\{\pm1\}.$ Hence, the same is true for $\mathcal{M}_{\lambda\lambda}^{\nu}(w_0).$ Let $\mathcal{M}_{\lambda\lambda}^{\nu}(w_0,\pm1)$ denote the eigenspaces of the operator $\widetilde{w_0}$ on $\mathcal{M}_{\lambda\lambda}^{\nu}(w_0)$ corresponding to the eigenvalues $\pm 1$ respectively.

Panyushev-Yakimova \cite{PY} prove the following theorem in this special case. 
\begin{theorem}
Let $\lie g$ be a finite-dimensional semisimple Lie algebra, $\lambda$ be a self-dual dominant weight and $\nu$ be a dominant weight. Let $\epsilon_\lambda$ be the Frobenius-Schur indicator of $V(\lambda)$. Then,
\[
c^{\nu}_{\lambda\lambda}(+)=\dim \mathcal{M}_{\lambda\lambda}^{\nu}(w_0,\epsilon_\lambda)
\qquad
c^{\nu}_{\lambda\lambda}(-)=\dim\mathcal{M}_{\lambda\lambda}^{\nu}(w_0,-\epsilon_\lambda).
\] \qed
\end{theorem}
\begin{rem}
    Note the "twist" by $\epsilon_\lambda$ in the expression of the multiplicities. 
\end{rem}

\subsection{}
Our main theorem generalises the result due to Panyushev and Yakimova in several directions.
\begin{itemize}
    \item[1.] We allow $\lie g$ to be finite dimensional semisimple or \textit{symmetric Kac--Moody}.
    \item[2.] We allow $\lambda$ to be any dominant weight. That is, we allow arbitrary $V(\lambda)\in \catoint$ which is \textit{not necessarily self-dual} (in the finite dimensional case).
    \item[3.] We extend the result to Kostant--Kumar modules and \textit{any involution $w\in W$.} Letting $w=w_0$ recovers the result for the full tensor product in the finite dimensional case.
\end{itemize}

Let $w\in W$ be an involution. Consider the Kostant--Kumar module $K(\lambda,w,\lambda)=\mathfrak{Ug} (v_\lambda\otimes v_{w\lambda})$ where we may assume without loss of generality that $\widetilde{w}(v_\lambda)=v_{w\lambda}.$ By Corollary ~\ref{cor:eigen}, the eigenvalues of $\widetilde{w}$ on $V(\nu)_{\lambda+w\lambda}$ belong to $\{\pm1\}.$ Hence, the same is true for the $\widetilde{w}-$ invariant subspace $\mathcal{M}_{\lambda\lambda}^{\nu}(w)$ (see Proposition ~\ref{prop:invar}). Let the eigenspaces of $\widetilde{w}$ on $\mathcal{M}_{\lambda\lambda}^{\nu}(w)$ corresponding to eigenvalues $\pm1$ be denoted by $\mathcal{M}_{\lambda\lambda}^{\nu}(w,\pm1)$ respectively.

Recall from Proposition ~\ref{proptitslift} that $\widetilde{w}^2$ acts on the weight space $\mathbb{C}v_\lambda$ with the eigenvalue $(-1)^{\langle\lambda,\rho^\vee-w\rho^\vee\rangle},$ which is equal to $\varepsilon_{w\lambda}.$

The following is the main theorem of this paper.

\begin{theorem}\label{main}
     Let $\mathfrak{g}$ be a finite dimensional semisimple Lie algebra or a symmetric Kac--Moody Lie algebra. Let $\lambda$ and $\nu$ be dominant weights and $w$ be an element of the Weyl group such that $w=w^{-1}$. Then, $$c_{\lambda\lambda}^{\nu}(w,+)=\dim \mathcal{M}_{\lambda\lambda}^{\nu}(w,\varepsilon_{w\lambda}) \quad c_{\lambda\lambda}^{\nu}(w,-)=\dim \mathcal{M}_{\lambda\lambda}^{\nu}(w,-\varepsilon_{w\lambda}).$$
\end{theorem}
\begin{proof}
   Let $\tau$ denote the swap map in $V(\lambda)\otimes V(\lambda)$ that leaves $K(\lambda,w,\lambda)$ invariant. $K^+(\lambda, w,\lambda)$ and $K^-(\lambda, w,\lambda)$ are the eigenspaces of $\tau$ inside $K(\lambda, w,\lambda)$ corresponding to eigenvalues $+1$ and $-1$ respectively. Hence, $ \text{Hom}_{\lie g}(K^+(\lambda,w,\lambda),V(\nu))$ can be identified with $\{f\in \text{Hom}_{\lie g}(K(\lambda,w,\lambda),V(\nu))\;:\: f\tau=f\}$ and
   $f\tau=f$ if and only if  $$f(v_\lambda\otimes v_{w\lambda})=f(v_{w\lambda}\otimes v_\lambda).$$ Also, since $f$ is $\lie g$ linear, $fx=xf$ for all $x\in \text{exp}\;\lie g.$ Now, $f(v_\lambda\otimes v_{w\lambda})\in \mathcal{M}_{\lambda\lambda}^{\nu}(w)$ and $ (v_{w\lambda}\otimes v_\lambda)= \varepsilon_{w\lambda}\;\widetilde{w} (v_\lambda\otimes v_{w\lambda}). $ Therefore, $f(v_\lambda\otimes v_{w\lambda})=\varepsilon_{w\lambda}\;\widetilde{w}f(v_\lambda\otimes v_{w\lambda}).$ Since $\varepsilon_{w\lambda} $ has value either $1$ or $-1$, we have $\varepsilon_{w\lambda} f(v_\lambda\otimes v_{w\lambda})=\widetilde{w}f(v_\lambda\otimes v_{w\lambda}),$ which implies $f(v_\lambda\otimes v_{w\lambda})$ is an eigenvector for the action of $\widetilde{w}$ corresponding to the eigenvalue $\varepsilon_{w\lambda}.$ Conversely, given an eigenvector $v'$ in $\mathcal{M}_{\lambda\lambda}^{\nu}$ for eigenvalue $\varepsilon_{w\lambda},$ we can define $f$ such that $f(v_\lambda\otimes v_{w\lambda})= v'$ and $f\tau=f.$ \\
   Similarly, $V(\nu)$ is an irreducible in the decomposition of $K^-(\lambda, w,\lambda)$ if and only if $f(v_\lambda\otimes v_{w\lambda})\in \mathcal{M}_{\lambda\lambda}^\nu(w)$ is an eigenvector for the action of $\widetilde{w}$ corresponding to the eigenvalue $-\varepsilon_{w\lambda}.$
\end{proof}

\begin{example}
Consider the Lie algebra $\widehat{\lie{sl}_2}.$ Define $w\in W$ with reduced expression $s_0s_1s_0s_1\ldots s_0$ of odd length $m$ and let $\Lambda_0$ be the fundamental weight satisfying $\langle\Lambda_0, \alpha_0^\vee\rangle=1$ and $\langle\Lambda_0, \alpha_1^\vee\rangle=0.$ Then, the multiplicity of an irreducible $V(\nu)$ in $K^+(\Lambda_0,w,\Lambda_0)$ is given by the dimension of the space 
    $$\bigl \{v\in V(\nu)_{2\Lambda_0-(m+1)\alpha_0/2-(m^2-1)\delta/4}:\widetilde{w}v=(-1)^{(m+1)/2}v,$$
    $$\; e_{\alpha_0+k\delta}^{\text{max}(1,m+1-k)}v = 0, e_{\alpha_1+l\delta}v = 0,\;k\ge 0, l\ge0\; \text{and}\; \lie g_\beta v=0, \text{for}\; \beta \in \Delta_{im}^+\bigr\}.$$
\end{example}

\section{Combinatorial interpretation}
In type $A$, Carr{\'e} and Leclerc \cite{CL}, gave an expression in terms of domino tableaux, of the multiplicity of an irreducible in the decomposition of $\text{Sym}^2V(\lambda)$ and $\text{Alt}^2V(\lambda)$. 

Combining Corollary 5.5 of \cite{CL} and our theorem ~\ref{main} gives the following combinatorial interpretation of the dimensions of the split multiplicity spaces.

\begin{theorem}
Consider the Lie algebra $\lie g=\lie{sl}_{n+1}$ and let $\lambda$ and $\nu$ be dominant weights. Then,
\begin{itemize}

   \item[a)] $\dim \{v \in V(\nu)_{\lambda+w_0\lambda}: \widetilde{w_0}v= (-1)^{\langle \lambda, 2\rho^\vee\rangle}v,  e_{i}^{-\langle w_0\lambda, \alpha_i^\vee \rangle + 1} v =0 \text{ for all simple roots } \alpha_i \}$ is equal to the number of Yamanouchi domino tableaux of shape \\ $(2\lambda_1,2\lambda_1, 2\lambda_2,2\lambda_2,\ldots,2\lambda_n,2\lambda_n)$ and evaluation $\nu$ whose number of horizontal dominoes is $\equiv0\mod 4.$
   \item [b)] $\dim \{v \in V(\nu)_{\lambda+w_0\lambda}: \widetilde{w_0}v= (-1)^{\langle \lambda, 2\rho^\vee\rangle+1}v,  e_{i}^{-\langle w_0\lambda, \alpha_i^\vee \rangle + 1} v =0 \text{ for all simple roots } \alpha_i\}$ is equal to the number of Yamanouchi domino tableaux of shape \\ $(2\lambda_1,2\lambda_1, 2\lambda_2,2\lambda_2,\ldots,2\lambda_n,2\lambda_n)$ and evaluation $\nu$ whose number of horizontal dominoes is $\equiv2\mod 4.$
   \end{itemize}
\end{theorem}

\section{Acknowledgements}
The author thanks Sankaran Viswanath for suggesting the problem and for numerous helpful discussions. 

\printbibliography

@article {jeralds-kumar,
    AUTHOR = {Jeralds, Samuel and Kumar, Shrawan},
     TITLE = {Root components for tensor product of affine {K}ac-{M}oody
              {L}ie algebra modules},
   JOURNAL = {Represent. Theory},
  FJOURNAL = {Representation Theory. An Electronic Journal of the American
              Mathematical Society},
    VOLUME = {26},
      YEAR = {2022},
     PAGES = {825--858},
      ISSN = {1088-4165},
   MRCLASS = {17B67 (14M15 14M17 20G44 22E47)},
  MRNUMBER = {4457456},
MRREVIEWER = {Yongjie\ Wang},
       DOI = {10.1090/ert/617},
       URL = {https://doi.org/10.1090/ert/617},
}

@article {kostant-tensmult,
    AUTHOR = {Kostant, Bertram},
     TITLE = {Clifford algebra analogue of the {H}opf-{K}oszul-{S}amelson
              theorem, the {$\rho$}-decomposition {$C(\lie g)={\rm End}\,
              V_\rho \otimes C(P)$}, and the {$\lie g$}-module structure of
              {$\bigwedge \lie g$}},
   JOURNAL = {Adv. Math.},
  FJOURNAL = {Advances in Mathematics},
    VOLUME = {125},
      YEAR = {1997},
    NUMBER = {2},
     PAGES = {275--350},
      ISSN = {0001-8708,1090-2082},
   MRCLASS = {17B20 (22E60)},
  MRNUMBER = {1434113},
MRREVIEWER = {M.\ R.\ Sepanski},
       DOI = {10.1006/aima.1997.1608},
       URL = {https://doi.org/10.1006/aima.1997.1608},
}

@article{CL,
  author  = {Carr{\'e}, Christophe and Leclerc, Bernard},
  title   = {Splitting the square of a Schur function into its symmetric and antisymmetric parts},
  journal = {Journal of Algebraic Combinatorics},
  year    = {1995},
  volume  = {4},
  number  = {3},
  pages   = {201--231},
  doi     = {https://doi.org/10.1023/A%3A1022475927626}
}

@article{GE,
  author  = {Gould, M. D. and Edwards, S. A.},
  title   = {Enveloping algebra annihilators and projection techniques for
             finite dimensional cyclic modules of a semisimple {L}ie algebra},
  journal = {Journal of Mathematical Physics},
  volume  = {25},
  pages   = {2848--2855},
  year    = {1984},
  doi     = {10.1063/1.526033},
}

@book{Hum1,
  author    = {Humphreys, James E.},
  title     = {Introduction to {L}ie Algebras and Representation Theory},
  series    = {Graduate Texts in Mathematics},                
  volume    = {9},                                            
  publisher = {Springer},
  address   = {New York, NY},
  year      = {1973},
  doi       = {10.1007/978-1-4612-6398-2},
}

@book{Kac,
  author    = {Kac, V. G.},
  title     = {Infinite-Dimensional {L}ie Algebras},
  edition   = {3},
  publisher = {Cambridge University Press},
  address   = {Cambridge},
  year      = {1990},
  doi       = {10.1017/CBO9780511626234},
}

@article{KRV,
  author  = {Kushwaha, Mrigendra Singh and Raghavan, K. N. and
             Viswanath, Sankaran},
  title   = {A study of {K}ostant--{K}umar modules via {L}ittelmann paths},
  journal = {Advances in Mathematics},
  volume  = {381},
  pages   = {107614},
  year    = {2021},
  issn    = {0001-8708},
  doi     = {10.1016/j.aim.2021.107614},
}

@misc{KRV2,
      title={The saturation property for refined Littlewood-Richardson coefficients}, 
      author={Mrigendra Singh Kushwaha and K. N. Raghavan and Sankaran Viswanath},
      year={2025},
      eprint={2204.03399},
      archivePrefix={arXiv},
      primaryClass={math.RT},
      url={https://arxiv.org/abs/2204.03399}, 
}

@book{Kum3,
  author    = {Shrawan Kumar},
  title     = {Kac-Moody Groups, their Flag Varieties and Representation Theory},
  publisher = {Birkh{\"a}user Boston},
  series    = {Progress in Mathematics},
  volume    = {204},
  year      = {2002},
  isbn      = {0-8176-4227-7},
  doi       = {10.1007/978-1-4612-0105-2}
}

@misc{MRV,
      title={Kostant--Kumar modules: presentation and multiplicities}, 
      author={Manika Gupta and K. N. Raghavan and Sankaran Viswanath},
      year={2026},
      eprint={2608.14098},
      archivePrefix={arXiv},
      primaryClass={math.RT},
      url={https://arxiv.org/abs/2608.14098}, 
}

@article{PRV,
  author    = {K. R. Parthasarathy and R. Ranga Rao and V. S. Varadarajan},
  title     = {Representations of complex semi-simple {Lie} groups and {Lie} algebras},
  journal   = {Annals of Mathematics},
  series    = {2},
  volume    = {85},
  pages     = {383--429},
  year      = {1967},
  doi       = {10.2307/1970351}
}

@article{PY,
  author  = {Panyushev, D. I. and Yakimova, O. S.},
  title   = {The PRV-formula for tensor product decompositions and its applications},
  journal = {Functional Analysis and Its Applications},
  year    = {2008},
  volume  = {42},
  number  = {1},
  pages   = {45--52},
  doi     = {10.1007/s10688-008-0005-7}
}

@book{S1,
  author    = {Steinberg, Robert},
  title     = {Lectures on Chevalley Groups},
  series    = {University Lecture Series},
  volume    = {66},
  publisher = {American Mathematical Society},
  address   = {Providence, Rhode Island},
  year      = {2016},
  note      = {Notes prepared by John Faulkner and Robert Wilson; revised and corrected edition},
  isbn      = {978-1-4704-3105-1},
  url       = {https://www.ams.org/books/ulect/066/}
}

@article{T1,
  title={Normalisateurs de tores I. Groupes de coxeter {\'e}tendus},
  author={Jacques Tits},
  journal={Journal of Algebra},
  year={1966},
  volume={4},
  pages={96-116},
  url={https://api.semanticscholar.org/CorpusID:121195695}
}

\end{document}